\input amstex
\documentstyle{amsppt}
%
\catcode`@=11
\redefine\output@{%
  \def\break{\penalty-\@M}\let\par\endgraf
  \ifodd\pageno\global\hoffset=105pt\else\global\hoffset=8pt\fi  
  \shipout\vbox{%
    \ifplain@
      \let\makeheadline\relax \let\makefootline\relax
    \else
      \iffirstpage@ \global\firstpage@false
        \let\rightheadline\frheadline
        \let\leftheadline\flheadline
      \else
        \ifrunheads@ 
        \else \let\makeheadline\relax
        \fi
      \fi
    \fi
    \makeheadline \pagebody \makefootline}%
  \advancepageno \ifnum\outputpenalty>-\@MM\else\dosupereject\fi
}
\def\Beta{\mathchar"0\hexnumber@\rmfam 42}
\catcode`\@=\active
\nopagenumbers
\chardef\textvolna='176

\chardef\bigalpha='013
\def\negskp{\hskip -2pt}

\chardef\degree="5E
\def\compos{\,\raise 1pt\hbox{$\sssize\circ$} \,}

\def\blue#1{#1}

\catcode`#=11\def\diez{#}\catcode`#=6
\catcode`&=11\catcode`&=4
\catcode`_=11\def\podcherkivanie{_}\catcode`_=8
\catcode`\^=11\catcode`\^=7
\catcode`~=11\catcode`~=\active
\def\mycite#1{\cite{\blue{#1}}\immediate\special{ps:
     ShrHPSdict begin /ShrBORDERthickness 0 def}}
\def\myciterange#1#2#3#4{\cite{\blue{#2#3#4}}\immediate\special{ps:
     ShrHPSdict begin /ShrBORDERthickness 0 def}}
\def\mytag#1{%
    \tag#1}
\def\mythetag#1{\thetag{\blue{#1}}\immediate\special{ps:
     ShrHPSdict begin /ShrBORDERthickness 0 def}}
\def\myrefno#1{\no#1}
\def\myhref#1#2{\blue{#2}\immediate\special{ps:
     ShrHPSdict begin /ShrBORDERthickness 0 def}}
\def\myEarXivlink{\myhref{http://arXiv.org}{http:/\negskp/arXiv.org}}

\def\mytheorem#1{\csname proclaim\endcsname{Theorem #1}}
\def\mytheoremwithtitle#1#2{\csname proclaim\endcsname{Theorem #1#2}}
\def\mythetheorem#1{\blue{#1}\immediate\special{ps:
     ShrHPSdict begin /ShrBORDERthickness 0 def}}
\def\mylemma#1{\csname proclaim\endcsname{Lemma #1}}
\def\mylemmawithtitle#1#2{\csname proclaim\endcsname{Lemma #1#2}}

\def\mycorollary#1{\csname proclaim\endcsname{Corollary #1}}

\def\mydefinition#1{\definition{Definition #1}}
\def\mythedefinition#1{\blue{#1}\immediate\special{ps:
     ShrHPSdict begin /ShrBORDERthickness 0 def}}
\def\myconjecture#1{\csname proclaim\endcsname{Conjecture #1}}
\def\myconjecturewithtitle#1#2{\csname proclaim\endcsname{Conjecture #1#2}}

\def\myproblem#1{\csname proclaim\endcsname{Problem #1}}
\def\myproblemwithtitle#1#2{\csname proclaim\endcsname{Problem #1#2}}


\pagewidth{360pt}
\pageheight{606pt}
\topmatter
\title
On quartic forms associated with cubic transformations of the real plane.
\endtitle
\rightheadtext{On quartic forms associated with cubic transformations \dots}
\author
Ruslan Sharipov
\endauthor
\address Bashkir State University, 32 Zaki Validi street, 450074 Ufa, Russia
\endaddress
\email
\myhref{mailto:r-sharipov\@mail.ru}{r-sharipov\@mail.ru}
\endemail
\abstract
     A polynomial transformation of the real plane $\Bbb R^2$ is a mapping
$\Bbb R^2\to\Bbb R^2$ given by two polynomials of two variables. Such a
transformation is called cubic if the degrees of its polynomials are not
greater than three. It turns out that cubic transformations are associated 
with some binary and quaternary quartic forms. In the present paper these 
forms are defined and studied. 
\endabstract
\endtopmatter
\TagsOnRight
\document

\head
1. Introduction.
\endhead
     Polynomial mappings $f\!:\,\Bbb C^2\to\Bbb C^2$ are of interests from
many points of view (see \myciterange{1}{1}{--}{5}) most of which are associated
with the Jacobian conjecture (see \mycite{6}). Polynomial mappings $f\!:\,\Bbb R^2
\to\Bbb R^2$ are also considered (see \mycite{7}, \mycite{8}) in connection with
the real Jacobian conjecture. The rational real Jacobian conjecture is its 
generalization (see \mycite{9}). My interest is concentrated on cubic polynomial
transformations of $\Bbb R^2$ due to their application to the perfect cuboid 
problem (see \mycite{10} and \mycite{11}).\par 
     Let $f\!:\,\Bbb R^2\to\Bbb R^2$ be a cubic transformation of $\Bbb R^2$ 
and let $x=(x^1,x^2)$ be a point\footnotemark\ of $\Bbb R^2$. Assume that 
$y=(y^1,y^2)$ is the image of $x$ under this transformation. Then 
\footnotetext{\ We use upper and lower indices according to Einstein's 
tensorial notation (see \mycite{12})} 
$$
\hskip -2em
y^{\kern 0.5pt i}=\sum^2_{m=1}\sum^2_{n=1}\sum^2_{p=1}
F^{\kern 0.5pt i}_{m\kern 0.1pt np}\,x^m\,x^n\,x^{\kern 0.1pt p}
+\ldots.
\mytag{1.1}
$$
We use dots in order to designate lower degree terms, i\.\,e\. quadratic, 
linear and constant terms. The formula \mythetag{1.1} is a coordinate 
presentation of the transforma\-tion $f\!:\,\Bbb R^2\to\Bbb R^2$. The coefficients 
$F^i_{m\kern 0.1pt np}$ in \mythetag{1.1} are components of a tensor
of the type $(1,3)$ (see definition in \mycite{13}).\par
     Along with \mythetag{1.1}, we consider linear transformations 
$\varphi\!:\,\Bbb R^2\to\Bbb R^2$ of the form
$$
\hskip -2em
y^{\kern 0.5pt i}=\sum^2_{m=1}T^{\kern 0.5pt i}_m\,x^m+a^i.
\mytag{1.2}
$$
\mydefinition{1.1} Two cubic mappings $f\!:\,\Bbb R^2\to \Bbb R^2$ and
$\tilde f\!:\,\Bbb R^2\to \Bbb R^2$ are called equivalent if there are two 
invertible linear transformations $\varphi_1\!:\,\Bbb R^2\to \Bbb R^2$ and
$\varphi_2\!:\,\Bbb R^2\to \Bbb R^2$ such that \pagebreak $\varphi_1\compos f=
\tilde f\compos\varphi_2$.
\enddefinition
     Invertibility of a linear mapping $\varphi\!:\,\Bbb R^2\to\Bbb R^2$
means invertibility of its matrix $T$ in \mythetag{1.2}. Assume that the 
mappings $\varphi_1$ and $\varphi_2$ in Definition~\mythedefinition{1.1}
are given by the matrices $T[1]$ and $T[2]$ and assume that $S[1]=T[1]^{-1}$. 
Then the tensors $F$ and $\tilde F$ associated with the cubic transformations
$f$ and $\tilde f$ are related to each other as 
$$
\hskip -2em
F^{\kern 0.5pt i}_{m\kern 0.1pt np}=
\sum^2_{\tilde{\kern -2pt i}=1}\sum^2_{\tilde m=1}\sum^2_{\tilde n=1}
\sum^2_{\tilde p=1}
\tilde F^{\kern 2.5pt\tilde{\kern -2pt i}}_{\tilde m\kern 0.1pt\tilde n
\tilde p}\ S^{\kern 0.5pt i}_{\kern 1.5pt\tilde{\kern -2pt i}}[1]
\ T^{\kern 0.5pt\tilde m}_m[2]\ T^{\kern 0.5pt\tilde n}_n[2]
\ T^{\kern 0.5pt\tilde p}_p[2].
\mytag{1.3}
$$\par 
The formula \mythetag{1.3} is slightly different from the transformation formula 
for components of a tensor under a change of Cartesian coordinates in $\Bbb R^2$
(see \mycite{13}). The main goal of the present paper is to introduce some tensorial
invariants associated with $F$ and study their behavior under non-tensorial 
transformations given by the formula \mythetag{1.3}. These invariants are analogous 
to those considered in \mycite{14} in the case of quadratic transformations of 
the real plane $\Bbb R^2$.\par
\head
2. Determinants and their collective behavior. 
\endhead
     The tensor $F$ in \mythetag{1.1} is symmetric with respect to its lower
indices. Taking into account this symmetry we can write \mythetag{1.1} as 
$$
\hskip -2em
\aligned
&y^1=F^{\kern 1pt 1}_{111}\,(x^1)^3+3\,F^{\kern 1pt 1}_{112}\,(x^1)^2\,x^2
+3\,F^{\kern 1pt 1}_{122}\,x^1\,(x^2)^2+F^{\kern 1pt 1}_{222}\,(x^2)^3,\\
&y^2=F^{\kern 1pt 2}_{111}\,(x^1)^3+3\,F^{\kern 1pt 2}_{112}\,(x^1)^2\,x^2
+3\,F^{\kern 1pt 2}_{122}\,x^1\,(x^2)^2+F^{\kern 1pt 2}_{222}\,(x^2)^3.
\endaligned
\mytag{2.1}
$$
Like in \mycite{14}, using \mythetag{2.1}, we consider the following determinants 
$$
\xalignat 2
&\hskip -2em
G_{1111}=\vmatrix F^{\kern 1pt 1}_{111} & F^{\kern 1pt 1}_{112}\\
\vspace{2ex}
F^{\kern 1pt 2}_{111} & F^{\kern 1pt 2}_{112}\endvmatrix,
&&G_{1112}=\vmatrix F^{\kern 1pt 1}_{111} & F^{\kern 1pt 1}_{122}\\
\vspace{2ex}
F^{\kern 1pt 2}_{111} & F^{\kern 1pt 2}_{122}\endvmatrix,\\
\vspace{2ex}
&\hskip -2em
G_{1122}=\vmatrix F^{\kern 1pt 1}_{111} & F^{\kern 1pt 1}_{222}\\
\vspace{2ex}
F^{\kern 1pt 2}_{111} & F^{\kern 1pt 2}_{222}\endvmatrix,
&&G_{1212}=\vmatrix F^{\kern 1pt 1}_{112} & F^{\kern 1pt 1}_{122}\\
\vspace{2ex}
F^{\kern 1pt 2}_{112} & F^{\kern 1pt 2}_{122}\endvmatrix,
\mytag{2.2}\\
\vspace{2ex}
&\hskip -2em
G_{1222}=\vmatrix F^{\kern 1pt 1}_{112} & F^{\kern 1pt 1}_{222}\\
\vspace{2ex}
F^{\kern 1pt 2}_{112} & F^{\kern 1pt 2}_{222}\endvmatrix,
&&G_{2222}=\vmatrix F^{\kern 1pt 1}_{122} & F^{\kern 1pt 1}_{222}\\
\vspace{2ex}
F^{\kern 1pt 2}_{122} & F^{\kern 1pt 2}_{222}\endvmatrix.
\endxalignat
$$
Using \mythetag{2.2}, we define several quartic forms in $\Bbb R^2$ and
in $\Bbb R^4$. Three of them are 
$$
\gather
\hskip -2em
\gathered
\omega[1]=G_{1111}\,(z^1)^4+2\,G_{1112}\,(z^1)^3\,z^2\,+\\
+\,(3\,G_{1212}+G_{1122})\,(z^1)^2\,(z^2)^2+2\,G_{1222}\,z^1\,(z^2)^3
+G_{2222}\,(z^2)^4,
\endgathered\\
\vspace{2ex}
\hskip -2em
\gathered
\omega[2]=2\,G_{1111}\,(z^1)^3\,z^3+G_{1112}\,(z^1)^3\,z^4
+3\,G_{1112}\,(z^1)^2\,z^2\,z^3\,+\\
+\,(3\,G_{1212}+G_{1122})\,(z^1)^2\,z^2\,z^4
+(3\,G_{1212}+G_{1122})\,z^1\,(z^2)^2\,z^3\,+\,\\
+\,3\,G_{1222}\,z^1\,(z^2)^2\,z^4
+G_{2222}\,\,(z^2)^3\,z^3+2\,G_{2222}\,\,(z^2)^3\,z^4,
\endgathered
\mytag{2.3}\\
\vspace{2ex}
\hskip -2em
\gathered
\omega[3]=3\,G_{1111}\,(z^1)^2\,(z^3)^2+3\,G_{1112}\,(z^1)^2\,z^3\,z^4
+G_{1122}\,(z^1)^2\,(z^4)^2\,+\\
+\,3\,G_{1112}\,z^1\,z^2\,(z^3)^2
+(9\,G_{1212}+G_{1122})\,z^1\,z^2\,z^3\,z^4+3\,G_{1222}\,z^1\,z^2\,(z^4)^2\,+\\
+\,G_{1122}\,(z^2)^2\,(z^3)^2
+3\,G_{1222}\,\,(z^2)^2\,z^3\,z^4+3\,G_{2222}\,\,(z^2)^2\,(z^4)^2.
\endgathered
\endgather
$$
The other three quartic forms are
$$
\gather
\hskip -2em
\gathered
\omega[4]=G_{1111}\,(z^1)^2\,(z^3)^2+G_{1112}\,(z^1)^2\,z^3\,z^4
+G_{1212}\,(z^1)^2\,(z^4)^2\,+\\
+\,G_{1112}\,z^1\,z^2\,(z^3)^2
+\,(G_{1212}+G_{1122})\,z^1\,z^2\,z^3\,z^4
+G_{1222}\,z^1\,z^2\,(z^4)^2\,+\\
+\,G_{1212}\,(z^2)^2\,(z^3)^2+G_{1222}\,(z^2)^2\,z^3\,z^4 
+G_{1222}\,(z^2)^2\,(z^4)^2,
\endgathered\\
\vspace{2ex}
\hskip -2em
\gathered
\omega[5]=2\,G_{1111}\,z^1\,(z^3)^3+G_{1112}\,z^2\,(z^3)^3
+3\,G_{1112}\,z^1\,(z^3)^2\,z^4\,+\\
+\,(3\,G_{1212}+G_{1122})\,z^2\,(z^3)^2\,z^4
+(3\,G_{1212}+G_{1122})\,z^1\,z^3\,(z^4)^2\,+\\
+\,3\,G_{1222}\,z^2\,z^3\,(z^4)^2
+G_{1222}\,z^1\,(z^4)^3
+2\,G_{2222}\,z^2\,(z^4)^3,
\endgathered
\mytag{2.4}\\
\vspace{2ex}
\hskip -2em
\gathered
\omega[6]=G_{1111}\,(z^3)^4+2\,G_{1112}\,(z^3)^3\,z^4\,+\\
+\,(3\,G_{1212}+G_{1122})\,(z^3)^2\,(z^4)^2
+2\,G_{1222}\,z^3\,(z^4)^3+G_{2222}\,(z^4)^4.
\endgathered
\endgather
$$
The quantities $z^1$, $z^2$, $z^3$, $z^4$ are interpreted as
components of two vectors 
$$
\xalignat 2
&\hskip -2em
\Vmatrix z^1\\ \vspace{1ex} z^2\endVmatrix\in\Bbb R^2,
&&\Vmatrix z^3\\ \vspace{1ex} z^4\endVmatrix\in\Bbb R^2
\mytag{2.5}
\endxalignat
$$
in $\omega[1]$ and $\omega[6]$ or as components of a single vector 
$$
\hskip -2em
\Vmatrix z^1\\ \vspace{1ex} z^2 \vspace{1ex} z^3\\ \vspace{1ex} z^4
\endVmatrix\in\Bbb R^2\oplus\Bbb R^2=\Bbb R^4
\mytag{2.6}
$$
in $\omega[2]$, $\omega[3]$, $\omega[4]$, and $\omega[5]$. Using the 
components of the vectors \mythetag{2.5}, we can write 
$\omega[1]$ and $\omega[6]$ from \mythetag{2.3} and \mythetag{2.4} in terms
of their components:
$$
\align
&\hskip -2em
\omega[1]=\sum^2_{i=1}\sum^2_{m=1}\sum^2_{n=1}\sum^2_{p=1}
\Omega[1]_{im\kern 0.1pt np}\,z^i\,z^m\,z^n\,z^p,
\mytag{2.7}\\
&\hskip -2em
\omega[6]=\sum^4_{i=3}\sum^4_{m=3}\sum^4_{n=3}\sum^4_{p=3}
\Omega[6]_{im\kern 0.1pt np}\,z^i\,z^m\,z^n\,z^p.
\mytag{2.8}
\endalign
$$
Similarly, using the components of the vector \mythetag{2.6}, we can write 
the quartic forms $\omega[2]$, $\omega[3]$, $\omega[4]$, and $\omega[5]$ in 
terms of their components:
$$
\hskip -2em
\aligned
&\omega[2]=\sum^4_{i=1}\sum^4_{m=1}\sum^4_{n=1}\sum^4_{p=1}
\Omega[2]_{im\kern 0.1pt np}\,z^i\,z^m\,z^n\,z^p,\\
&\omega[3]=\sum^4_{i=1}\sum^4_{m=1}\sum^4_{n=1}\sum^4_{p=1}
\Omega[3]_{im\kern 0.1pt np}\,z^i\,z^m\,z^n\,z^p,\\
&\omega[4]=\sum^4_{i=1}\sum^4_{m=1}\sum^4_{n=1}\sum^4_{p=1}
\Omega[4]_{im\kern 0.1pt np}\,z^i\,z^m\,z^n\,z^p,\\
&\omega[5]=\sum^4_{i=1}\sum^4_{m=1}\sum^4_{n=1}\sum^4_{p=1}
\Omega[5]_{im\kern 0.1pt np}\,z^i\,z^m\,z^n\,z^p.
\endaligned
\mytag{2.9}
$$\par
     Now assume that we perform a linear change of coordinates in 
$\Bbb R^2$. It is expressed by the following matrix formulas for the 
vectors \mythetag{2.5}:
$$
\xalignat 2
&\hskip -2em
\Vmatrix z^1\\ \vspace{1ex} z^2\endVmatrix=
\Vmatrix S^1_1 & S^1_2\\ \vspace{1ex} S^2_1 & S^2_2\endVmatrix\cdot
\Vmatrix\tilde z^1\\ \vspace{1ex}\tilde z^2\endVmatrix,
&&\Vmatrix z^3\\ \vspace{1ex} z^4\endVmatrix=
\Vmatrix S^1_1 & S^1_2\\ \vspace{1ex} S^2_1 & S^2_2\endVmatrix\cdot
\Vmatrix\tilde z^3\\ \vspace{1ex}\tilde z^4\endVmatrix.
\quad
\mytag{2.10}
\endxalignat
$$
For the vector \mythetag{2.6} the formulas \mythetag{2.10} imply
$$
\hskip -2em
\Vmatrix z^1\\ \vspace{1ex} z^2\\ \vspace{1ex} z^3\\ \vspace{1ex} z^4
\endVmatrix=\Vmatrix S^1_1 & S^1_2 & 0 & 0\\ 
\vspace{1ex} S^2_1 & S^2_2 & 0 & 0\\
\vspace{1ex} 0 & 0 & S^1_1 & S^1_2\\ 
\vspace{1ex} 0 & 0 & S^2_1 & S^2_2\endVmatrix\cdot \Vmatrix 
\tilde z^1\\ \vspace{1ex}\tilde  z^2\\ \vspace{1ex}\tilde z^3\\ 
\vspace{1ex}\tilde z^4\endVmatrix.
\mytag{2.11}
$$
The matrix $S$ used in \mythetag{2.10} is called a transition matrix
(see \mycite{12} or \mycite{13}). Let's denote through $\hat S$ the 
block-diagonal matrix in \mythetag{2.11}. It plays the role of a
transition matrix for the linear change of coordinates \mythetag{2.11}
in $\Bbb R^4=\Bbb R^2\oplus\Bbb R^2$.\par 
     Each linear change of coordinates implies some definite associated 
change of components for all tensors (see \mycite{13}). In the case of 
the tensor $F$ in \mythetag{1.1} we have
$$
\hskip -2em
\tilde F^{\kern 0.5pt i}_{m\kern 0.1pt np}=
\sum^2_{\tilde{\kern -2pt i}=1}\sum^2_{\tilde m=1}\sum^2_{\tilde n=1}
\sum^2_{\tilde p=1}
 F^{\kern 2.5pt\tilde{\kern -2pt i}}_{\tilde m\kern 0.1pt\tilde n
\tilde p}\ T^{\kern 0.5pt i}_{\kern 1.5pt\tilde{\kern -2pt i}}
\ S^{\kern 0.5pt\tilde m}_m\ S^{\kern 0.5pt\tilde n}_n
\ S^{\kern 0.5pt\tilde p}_p.
\mytag{2.12}
$$
Here $T=S^{-1}$. The formula \mythetag{2.12} has its inverse formula
$$
\hskip -2em
F^{\kern 0.5pt i}_{m\kern 0.1pt np}=
\sum^2_{\tilde{\kern -2pt i}=1}\sum^2_{\tilde m=1}\sum^2_{\tilde n=1}
\sum^2_{\tilde p=1}
\tilde F^{\kern 2.5pt\tilde{\kern -2pt i}}_{\tilde m\kern 0.1pt\tilde n
\tilde p}\ S^{\kern 0.5pt i}_{\kern 1.5pt\tilde{\kern -2pt i}}
\ T^{\kern 0.5pt\tilde m}_m\ T^{\kern 0.5pt\tilde n}_n
\ T^{\kern 0.5pt\tilde p}_p,
\mytag{2.13}
$$
which is very similar to the formula \mythetag{1.3}, though the meanings of 
the formulas \mythetag{1.3} and \mythetag{2.13} are quite different.\par
     Now we can use \mythetag{2.12} in \mythetag{2.2} instead of 
$F^{\kern 0.5pt i}_{m\kern 0.1pt np}$ and obtain six determinants 
$\tilde G_{1111}$, $\tilde G_{1112}$, $\tilde G_{1122}$, $\tilde G_{1212}$, 
$\tilde G_{1222}$, $\tilde G_{2222}$. Then we can use these determinants 
in \mythetag{2.3} and \mythetag{2.4} instead of $G_{1111}$, $G_{1112}$, 
$G_{1122}$, $G_{1212}$, $G_{1222}$, $G_{2222}$ simultaneously replacing
$z^1$, $z^2$, $z^3$, $z^4$ by $\tilde z^1$, $\tilde z^2$, $\tilde z^3$, 
$\tilde z^4$. As a result we get some expressions for $\omega[1]$, 
$\omega[2]$, $\omega[3]$, $\omega[4]$, $\omega[5]$, $\omega[6]$ through 
$F$, $T$, $S$, and $\tilde z$. On the other hand we can apply
\mythetag{2.10} or \mythetag{2.11} directly to \mythetag{2.3} and 
\mythetag{2.4}. It turns out that the results of these two ways of 
expressing $\omega[1]$, $\omega[2]$, $\omega[3]$, $\omega[4]$, $\omega[5]$, 
$\omega[6]$ through $F$, $T$, $S$, and $\tilde z$ do always coincide. We
write this fact as 
$$
\hskip -2em
\omega[q](\tilde F(F),\tilde z)=\omega[q](F,z(\tilde z)),
\ q=1,\,\ldots,\,6.
\mytag{2.14}
$$
We can also express this fact using the component notations \mythetag{2.7}, 
\mythetag{2.8} and \mythetag{2.9}:
$$
\allowdisplaybreaks
\gather
\hskip -2em
\tilde\Omega[1]_{\kern 1.5pt\tilde{\kern -2pt i}\tilde m\kern 0.1pt
\tilde n\tilde p}=\sum^2_{i=1}\sum^2_{m=1}\sum^2_{n=1}\sum^2_{p=1}
\Omega[1]_{{im\kern 0.1pt np}}\
S^{\kern 0.5pt i}_{\kern 1.5pt\tilde{\kern -2pt i}}
\ S^{\kern 0.5pt m}_{\tilde m}\ S^{\kern 0.5pt n}_{\tilde n}
\ S^{\kern 0.5pt p}_{\tilde p},
\mytag{2.15}\\
\hskip -2em
\tilde\Omega[6]_{\kern 1.5pt\tilde{\kern -2pt i}\tilde m\kern 0.1pt
\tilde n\tilde p}=\sum^4_{i=3}\sum^4_{m=3}\sum^4_{n=3}\sum^4_{p=3}
\Omega[6]_{{im\kern 0.1pt np}}\
S^{\kern 0.5pt i}_{\kern 1.5pt\tilde{\kern -2pt i}}
\ S^{\kern 0.5pt m}_{\tilde m}\ S^{\kern 0.5pt n}_{\tilde n}
\ S^{\kern 0.5pt p}_{\tilde p},
\mytag{2.16}\\
\displaybreak
\hskip -2em
\tilde\Omega[q]_{\kern 1.5pt\tilde{\kern -2pt i}\tilde m\kern 0.1pt
\tilde n\tilde p}=\sum^4_{i=1}\sum^4_{m=1}\sum^4_{n=1}\sum^4_{p=1}
\Omega[q]_{{im\kern 0.1pt np}}\
\hat S^{\kern 0.5pt i}_{\kern 1.5pt\tilde{\kern -2pt i}}
\ \hat S^{\kern 0.5pt m}_{\tilde m}\ \hat S^{\kern 0.5pt n}_{\tilde n}
\ \hat S^{\kern 0.5pt p}_{\tilde p},\ q=2,\,\ldots,\,5.
\mytag{2.17}
\endgather
$$
The formula \mythetag{2.14}, as well as the formulas \mythetag{2.15},
\mythetag{2.16}, and \mythetag{2.17}, means that the components of
the quartic forms $\omega[1]$, $\omega[2]$, $\omega[3]$, $\omega[4]$, 
$\omega[5]$, $\omega[6]$ exhibit true tensorial behavior under the linear 
change of coordinates \mythetag{2.10} and \mythetag{2.11}.\par
\head
3. Explicit formulas for components of quartic forms. 
\endhead
     In our case $\dim\Bbb R^2=2$. There is a fundamental pseudotensor $\bold d$ of 
the type $(0,2)$ and of the weight $-1$ in each two-dimensional linear vector
space $V$. Its components are given by the following skew-symmetric matrix 
$$
\hskip -2em
d_{ij}=\Vmatrix\format \r&\quad\l\\ 0 & 1\\
\vspace{1ex}-1 & 0\endVmatrix.
\mytag{3.1}
$$
The dual object for $\bold d$ is given by the same matrix \mythetag{3.1}:
$$
\hskip -2em
d^{\kern 1pt ij}=\Vmatrix\format \r&\quad\l\\ 0 & 1\\
\vspace{1ex}-1 & 0\endVmatrix.
\mytag{3.2}
$$
This dual object is denoted by the same symbol $\bold d$ as the initial one. 
It is a pseudotensor of the type $(2,0)$, its weight is equal to $1$. The
following definition is provided for reference purposes only. 
\mydefinition{3.1} A pseudotensor of the type $(r,s)$ and 
of the weight $m$ in $\Bbb R^2$ is a geometrical and/or physical object presented 
by an array of quantities $F^{i_1\ldots\,i_r}_{j_1\ldots\,j_s}$ transformed as follows 
under any linear change of coordinates like \mythetag{2.10}: 
$$
\hskip -2em
F^{i_1\ldots\,i_r}_{j_1\ldots\,j_s}=
(\det T)^m\sum\Sb p_1\ldots p_r\\ q_1\ldots q_s\endSb
S^{i_1}_{p_1}\ldots\,S^{i_r}_{p_r}\,\,
T^{q_1}_{j_1}\ldots\,T^{q_s}_{j_s}\,\,
\tilde F^{p_1\ldots\,p_r}_{q_1\ldots\,q_s}.
\mytag{3.3}
$$
\enddefinition
The definition~\mythedefinition{3.1} can be easily modified for the case
of the space $\Bbb R^4=\Bbb R^2\oplus\Bbb R^2$ using the matrices $\hat S$
and $\hat T=\hat S^{-1}$ instead of $S$ and $T=S^{-1}$. 
\mydefinition{3.2} A pseudotensor of the type $(r,s)$ and the weight $m$ 
in the space $\Bbb R^4=\Bbb R^2\oplus\Bbb R^2$ is a geometrical and/or physical 
object presented by an array $F^{i_1\ldots\,i_r}_{j_1\ldots\,j_s}$ transformed 
as follows under any linear change of coordinates like \mythetag{2.11}: 
$$
\hskip -2em
F^{i_1\ldots\,i_r}_{j_1\ldots\,j_s}=
(\det T)^m\sum\Sb p_1\ldots p_r\\ q_1\ldots q_s\endSb
\hat S^{i_1}_{p_1}\ldots\,\hat S^{i_r}_{p_r}\,\,
\hat T^{q_1}_{j_1}\ldots\,\hat T^{q_s}_{j_s}\,\,
\tilde F^{p_1\ldots\,p_r}_{q_1\ldots\,q_s}.
\mytag{3.4}
$$
\enddefinition
     Pseudotensors of the weight $m=0$ are known as tensors (see \mycite{13}),
e\.\,g\. the formula \mythetag{2.13} is a particular instance of the formula
\mythetag{3.3}. The formulas \mythetag{2.15} and \mythetag{2.16} can be 
transformed to special instances of the formula \mythetag{3.3}, while the
formula \mythetag{2.17} can be transformed to a special instance of the formula 
\mythetag{3.4}.\par
     The weights of the pseudotensors \mythetag{3.1} and \mythetag{3.2} are
opposite to each other. Therefore they can be used in order to build a
tensor. \pagebreak Using these two pseudoten\-sors and the tensor $F$ from 
\mythetag{1.1}, we define a new tensor with the components 
$$
\hskip -2em
\Omega_{im\kern 0.1pt np}=\frac{1}{2}\sum^2_{r_1=1}\sum^2_{r_2=1}
\sum^2_{s_1=1}\sum^2_{s_2=1}F^{\kern 0.5pt r_1}_{s_1im}
\,F^{\kern 0.5pt r_2}_{s_2np}\,d^{\kern 0.5pt s_1s_2}
\,d_{\kern 0.5pt r_1r_2}.
\mytag{3.5}
$$
\mytheorem{3.1} The components of the quartic form $\omega[1]$ in
\mythetag{2.7} are produced from the components of the tensor 
\mythetag{3.5} by symmetrizing them with respect to $i$, $m$, $n$,
$p$:
$$
\hskip -2em
\Omega[1]_{im\kern 0.1pt np}=\frac{1}{3}
\bigl(\Omega_{im\kern 0.1pt np}+\Omega_{m\kern 0.1pt n\kern 0.1pt ip}
+\Omega_{n\kern 0.1pt imp}\bigr).
\mytag{3.6}
$$
\endproclaim 
    The proof of Theorem~\mythetheorem{3.1} is pure computations. 
In particular, one can easily verify that the right hand side of
the formula \mythetag{3.6} is fully symmetric with respect to the
indices $i$, $m$, $n$, $p$.\par
    The form $\omega[6]$ in \mythetag{2.4} does coincide with the 
form $\omega[1]$ in \mythetag{2.3} up to the substitution of $z^3$ for
$z^1$ and $z^4$ for $z^2$. Its components in \mythetag{2.8} are given
by the formula 
$$
\Omega[6]_{im\kern 0.1pt np}=\Omega[1]_{(i-2)\,(m-2)\,(n-2)\,(p-2)}.
$$\par
    The forms $\omega[2]$, $\omega[3]$, $\omega[4]$, $\omega[5]$ are more
complicated. In order to serve $\omega[2]$ we need to extend the tensor 
\mythetag{3.5} from $\Bbb R^2$ to $\Bbb R^4=\Bbb R^2\oplus\Bbb R^2$. We do 
it as follows:
$$
\hskip -2em
\hat\Omega_{im\kern 0.1pt np}=
\cases \Omega_{im\kern 0.1pt np} & \text{if \ }i\leqslant2,
\ m\leqslant2,\ n\leqslant2,\ p\leqslant2,\\
\quad 0 & \text{in all other cases.}
\endcases 
\mytag{3.7}
$$
Apart from \mythetag{3.7} we define the exchange operator $\varepsilon$
given by the matrix
$$
\hskip -2em
\varepsilon^{\kern 0.5pt i}_j
=\Vmatrix
0 & 0 & 1 & 0\\
\vspace{0.3ex}
0 & 0 & 0 & 1\\
\vspace{0.3ex}
1 & 0 & 0 & 0\\
\vspace{0.3ex}
0 & 1 & 0 & 0
\endVmatrix.
\mytag{3.8}
$$
Then we combine the tensors \mythetag{3.7} and \mythetag{3.8} in the 
following way:
$$
\hskip -2em
\Omega[A]_{im\kern 0.1pt np}=\sum^4_{\tilde p=1}
2\,\hat\Omega_{im\kern 0.1pt n\tilde p}
\ \varepsilon^{\kern 0.5pt\tilde p}_p\,.
\mytag{3.9}
$$
\mytheorem{3.2} The components of the form $\omega[2]$ in
\mythetag{2.9} are produced from the components of the tensor 
\mythetag{3.9} by symmetrizing them with respect to $i$, $m$, $n$,
$p$:
$$
\gathered
\Omega[2]_{im\kern 0.1pt np}=\frac{1}{12}
\bigl(\Omega[A]_{i\kern 0.1pt m\kern 0.1pt np}
+\Omega[A]_{ip\kern 0.4pt m\kern 0.1pt n}
+\Omega[A]_{i\kern 0.1pt np\kern 0.4pt m}
+\Omega[A]_{i\kern 0.1pt n\kern 0.1pt mp}\,+\\
+\,\Omega[A]_{ip\kern 0.4pt n\kern 0.1pt m}
+\Omega[A]_{i\kern 0.1pt mp\kern 0.4pt n}
+\Omega[A]_{mp\kern 0.4pt i\kern 0.1pt n}
+\Omega[A]_{mnp\kern 0.4pt i}
+\Omega[A]_{p\kern 0.4pt ni\kern 0.1pt m}\,+\\
+\,\Omega[A]_{p\kern 0.4pt m\kern 0.1pt ni}
+\Omega[A]_{n\kern 0.1pt mip}
+\Omega[A]_{np\kern 0.4pt mi}\bigr).
\endgathered
$$
\endproclaim 
     In order to serve the forms $\omega[3]$ and $\omega[4]$
we define the following two tensors: 
$$
\allowdisplaybreaks
\gather
\hskip -2em
\Omega[B]_{im\kern 0.1pt np}=\sum^4_{\tilde m=1}\sum^4_{\tilde p=1}
\hat\Omega_{i\tilde m\kern 0.1pt n\tilde p}
\ \varepsilon^{\kern 0.5pt\tilde m}_m
\ \varepsilon^{\kern 0.5pt\tilde p}_p
+\sum^4_{\tilde n=1}\sum^4_{\tilde p=1}
2\,\hat\Omega_{im\kern 0.1pt\tilde n\tilde p}
\ \varepsilon^{\kern 0.5pt\tilde n}_n
\ \varepsilon^{\kern 0.5pt\tilde p}_p\,,
\mytag{3.10}\\
\displaybreak
\hskip -2em
\Omega[C]_{im\kern 0.1pt np}=\sum^4_{\tilde m=1}\sum^4_{\tilde p=1}
\hat\Omega_{i\tilde m\kern 0.1pt n\tilde p}
\ \varepsilon^{\kern 0.5pt\tilde m}_m
\ \varepsilon^{\kern 0.5pt\tilde p}_p\,.
\mytag{3.11}\\
\endgather
$$
\mytheorem{3.3} The components of the form $\omega[3]$ in
\mythetag{2.9} are produced from the components of the tensor 
\mythetag{3.10} by symmetrizing them with respect to $i$, $m$, $n$,
$p$:
$$
\gathered
\Omega[3]_{im\kern 0.1pt np}=\frac{1}{24}
\bigl(\Omega[B]_{i\kern 0.1pt m\kern 0.1pt np}
+\Omega[B]_{ip\kern 0.4pt m\kern 0.1pt n}
+\Omega[B]_{i\kern 0.1pt np\kern 0.4pt m}
+\Omega[B]_{i\kern 0.1pt n\kern 0.1pt mp}\,+\\
+\,\Omega[B]_{ip\kern 0.4pt n\kern 0.1pt m}
+\Omega[B]_{i\kern 0.1pt mp\kern 0.4pt n}
+\Omega[B]_{m\kern 0.1pt i\kern 0.1pt np}
+\Omega[B]_{mp\kern 0.4pt i\kern 0.1pt n}
+\Omega[B]_{m\kern 0.1pt np\kern 0.4pt i}\,+\\
+\,\Omega[B]_{m\kern 0.1pt n\kern 0.1pt ip}
+\Omega[B]_{mp\kern 0.4pt n\kern 0.1pt i}
+\Omega[B]_{m\kern 0.1pt ip\kern 0.4pt n}
+\Omega[B]_{p\kern 0.4pt i\kern 0.1pt m\kern 0.1pt n}
+\Omega[B]_{p\kern 0.4pt ni\kern 0.1pt m}\,+\\
+\,\Omega[B]_{p\kern 0.4pt m\kern 0.1pt ni}
+\Omega[B]_{p\kern 0.4pt m\kern 0.1pt i\kern 0.1pt n}
+\Omega[B]_{p\kern 0.4pt n\kern 0.1pt m\kern 0.1pt i}
+\Omega[B]_{p\kern 0.4pt i\kern 0.1pt n\kern 0.1pt m}
+\Omega[B]_{n\kern 0.1pt ip\kern 0.4pt m}\,+\\
+\,\Omega[B]_{n\kern 0.1pt m\kern 0.1pt ip}
+\Omega[B]_{np\kern 0.4pt m\kern 0.1pt i}
+\Omega[B]_{n\kern 0.1pt p\kern 0.4pt i\kern 0.1pt m}
+\Omega[B]_{n\kern 0.1pt m\kern 0.1pt p\kern 0.4pt i}
+\Omega[B]_{ni\kern 0.1pt m\kern 0.1pt p}\bigr).
\endgathered
$$
\endproclaim 
\mytheorem{3.4} The components of the form $\omega[4]$ in
\mythetag{2.9} are produced from the components of the tensor 
\mythetag{3.11} by symmetrizing them with respect to $i$, $m$, $n$,
$p$:
$$
\gathered
\Omega[4]_{im\kern 0.1pt np}=\frac{1}{12}
\bigl(\Omega[C]_{i\kern 0.1pt m\kern 0.1pt np}
+\Omega[C]_{ip\kern 0.4pt m\kern 0.1pt n}
+\Omega[C]_{i\kern 0.1pt np\kern 0.4pt m}
+\Omega[C]_{i\kern 0.1pt n\kern 0.1pt mp}\,+\\
+\,\Omega[C]_{ip\kern 0.4pt n\kern 0.1pt m}
+\Omega[C]_{i\kern 0.1pt mp\kern 0.4pt n}
+\Omega[C]_{m\kern 0.1pt i\kern 0.1pt np}
+\Omega[C]_{m\kern 0.1pt np\kern 0.4pt i}
+\Omega[C]_{mp\kern 0.4pt n\kern 0.1pt i}\,+\\
+\,\Omega[C]_{m\kern 0.1pt ip\kern 0.4pt n}
+\Omega[C]_{p\kern 0.4pt m\kern 0.1pt ni}
+\Omega[C]_{p\kern 0.4pt i\kern 0.1pt n\kern 0.1pt m}
\bigr).
\endgathered
$$
\endproclaim 
     The form $\omega[5]$ is similar to $\omega[2]$. It differs
from $\omega[2]$ by exchanging $z^1\longleftrightarrow z^3$ and 
$z^2\longleftrightarrow z^4$. Such an exchange is performed by 
means of the operator \mythetag{3.8}. Therefore, looking at 
\mythetag{3.9}, by analogy we define the tensor 
$$
\hskip -2em
\Omega[D]_{im\kern 0.1pt np}=\sum^4_{\tilde m=1}\sum^4_{\tilde n=1}
\sum^4_{\tilde p=1}
2\,\hat\Omega_{i\tilde m\kern 0.1pt\tilde  n\tilde p}
\ \varepsilon^{\kern 0.5pt\tilde m}_m
\ \varepsilon^{\kern 0.5pt\tilde n}_n
\ \varepsilon^{\kern 0.5pt\tilde p}_p\,.
\mytag{3.12}
$$
\mytheorem{3.5} The components of the form $\omega[5]$ in
\mythetag{2.9} are produced from the components of the tensor 
\mythetag{3.12} by symmetrizing them with respect to $i$, $m$, 
$n$, $p$:
$$
\gathered
\Omega[5]_{im\kern 0.1pt np}=\frac{1}{12}
\bigl(\Omega[D]_{i\kern 0.1pt m\kern 0.1pt np}
+\Omega[D]_{ni\kern 0.1pt m\kern 0.1pt p}
+\Omega[D]_{m\kern 0.1pt ni\kern 0.1pt p}
+\Omega[D]_{i\kern 0.1pt n\kern 0.1pt mp}\\
+\Omega[D]_{m\kern 0.1pt i\kern 0.1pt n\kern 0.1pt p}
+\Omega[D]_{n\kern 0.1pt m\kern 0.1pt i\kern 0.1pt p}
+\Omega[D]_{p\kern 0.4pt i\kern 0.1pt m\kern 0.1pt n}
+\Omega[D]_{m\kern 0.1pt p\kern 0.4pt i\kern 0.1pt n}
+\Omega[D]_{p\kern 0.4pt ni\kern 0.1pt m}\,+\\
+\,\Omega[D]_{i\kern 0.1pt p\kern 0.4pt n\kern 0.1pt m}
+\Omega[D]_{pm\kern 0.1pt n\kern 0.1pt i}
+\Omega[D]_{n\kern 0.1pt p\kern 0.4pt m\kern 0.1pt i}\bigr).
\endgathered
$$
\endproclaim 
Theorems~\mythetheorem{3.2}, \mythetheorem{3.3}, \mythetheorem{3.4},
\mythetheorem{3.5} are similar to each other. All of them are proved by
means of direct calculations.\par
\head
4. Behavior under left and right compositions\\
with linear transformations. 
\endhead
    The tensorial behavior of the quartic forms $\omega[1]$, $\omega[2]$, 
$\omega[3]$, $\omega[4]$, $\omega[5]$, $\omega[6]$ revealed in Section~2 is 
essential for understanding their nature. However, it is inessential for 
their prospective applications. In the present section we study their behavior
under left and right compositions with linear transformations of the form
\mythetag{1.2} introduced in \pagebreak Definition~\mythedefinition{1.1}.\par
     Assume that we removed the right composition in 
Definition~\mythedefinition{1.1} so that we have the left composition only, 
i\.\,e\. assume that two cubic mappings $f\!:\,\Bbb R^2\to \Bbb R^2$ and
$\tilde f\!:\,\Bbb R^2\to \Bbb R^2$ of the form \mythetag{1.1} are 
related to each other as 
$$
\hskip -2em
f=\varphi^{-1}\compos \tilde f,
\mytag{4.1}
$$
where $\varphi\!:\,\Bbb R^2\to\Bbb R^2$ is a linear transformation of the 
form \mythetag{1.2}. In the case of \mythetag{4.1} the formula \mythetag{1.3}
reduces to the following one:
$$
\hskip -2em
F^{\kern 0.5pt i}_{m\kern 0.1pt np}=
\sum^2_{\tilde{\kern -2pt i}=1}
\tilde F^{\kern 2.5pt\tilde{\kern -2pt i}}_{m\kern 0.1ptn
p}\ S^{\kern 0.5pt i}_{\kern 1.5pt\tilde{\kern -2pt i}}.
\mytag{4.2}
$$ 
Here $S=T^{-1}$ and $T$ is the matrix from \mythetag{1.2}. Applying 
\mythetag{4.2} to \mythetag{2.2} we derive 
$$
\xalignat 2
&\hskip -2em
G_{1111}=\det S\cdot\tilde G_{1111}, 
&&G_{1112}=\det S\cdot\tilde G_{1112},\\
&\hskip -2em
G_{1122}=\det S\cdot\tilde G_{1122}, 
&&G_{1212}=\det S\cdot\tilde G_{1212},
\mytag{4.3}\\
&\hskip -2em
G_{1222}=\det S\cdot\tilde G_{1222}, 
&&G_{2222}=\det S\cdot\tilde G_{2222}.
\endxalignat
$$
Then, applying \mythetag{4.3} to the components of the quartic forms
$\omega[1]$, $\omega[2]$, $\omega[3]$, $\omega[4]$, $\omega[5]$, $\omega[6]$
in \mythetag{2.7}, \mythetag{2.8}, and \mythetag{2.9}, we derive
$$
\hskip -2em
\Omega[q]_{im\kern 0.1pt np}
=\det S\cdot\tilde \Omega[q]_{im\kern 0.1pt np} 
\text{, \ where \ }q=1,],\ldots,\,6.
\mytag{4.4}
$$
This result is summarized in the following theorem.
\mytheorem{4.1} Under the left composition \mythetag{4.1} of a cubic 
transformation $\tilde f$ with the inverse of a liner transformation 
$\varphi$ in \mythetag{1.2} its associated quartic forms $\omega[1]$, 
$\omega[2]$, $\omega[3]$, $\omega[4]$, $\omega[5]$, $\omega[6]$ are 
transformed according to the formulas \mythetag{4.4}.
\endproclaim 
     Now assume that we removed the left composition in 
Definition~\mythedefinition{1.1} so that we have the right composition only, 
i\.\,e\. assume that two cubic mappings $f\!:\,\Bbb R^2\to \Bbb R^2$ and
$\tilde f\!:\,\Bbb R^2\to \Bbb R^2$ of the form \mythetag{1.1} are 
related to each other as 
$$
\hskip -2em
f=\tilde f\compos\varphi,
\mytag{4.5}
$$
where $\varphi\!:\,\Bbb R^2\to\Bbb R^2$ is a linear transformation of the 
form \mythetag{1.2}. In the case of \mythetag{4.5} the formula \mythetag{1.3}
reduces to the following one:
$$
\hskip -2em
F^{\kern 0.5pt i}_{m\kern 0.1pt np}=\sum^2_{\tilde m=1}\sum^2_{\tilde n=1}
\sum^2_{\tilde p=1}\tilde F^{\kern 2.5pt i}_{\tilde m\kern 0.1pt\tilde n
\tilde p}
\ T^{\kern 0.5pt\tilde m}_m\ T^{\kern 0.5pt\tilde n}_n
\ T^{\kern 0.5pt\tilde p}_p.
\mytag{4.6}
$$
Applying \mythetag{4.6} to \mythetag{2.2} and then to \mythetag{2.3}, one can 
derive the formulas 
$$
\allowdisplaybreaks
\gather
\hskip -2em
\Omega[1]_{{im\kern 0.1pt np}}
=\det T\,\cdot
\sum^2_{\kern 1.5pt\tilde{\kern -2pt i}=1}
\sum^2_{\tilde m=1}\sum^2_{\tilde n=1}\sum^2_{\tilde p=1}
\tilde\Omega[1]_{\kern 2pt\tilde{\kern -2pt i}\kern 0.2pt\tilde m\kern 0.1pt
\tilde n\tilde p}
\ T^{\kern 2.5pt\tilde{\kern -2pt i}}_i
\ T^{\kern 0.5pt\tilde m}_m\ T^{\kern 0.5pt\tilde n}_n
\ T^{\kern 0.5pt \tilde p}_p,
\mytag{4.7}\\
\hskip -2em
\Omega[6]_{{im\kern 0.1pt np}}
=\det T\,\cdot
\sum^4_{\kern 1.5pt\tilde{\kern -2pt i}=3}
\sum^4_{\tilde m=3}\sum^4_{\tilde n=3}\sum^4_{\tilde p=3}
\tilde\Omega[6]_{\kern 2pt\tilde{\kern -2pt i}\kern 0.2pt\tilde m\kern 0.1pt
\tilde n\tilde p}
\ \hat T^{\kern 2.5pt\tilde{\kern -2pt i}}_i
\ \hat T^{\kern 0.5pt\tilde m}_m\ \hat T^{\kern 0.5pt\tilde n}_n
\ \hat T^{\kern 0.5pt \tilde p}_p,
\mytag{4.8}\\
\displaybreak
\hskip -2em
\Omega[q]_{{im\kern 0.1pt np}}
=\det T\,\cdot
\sum^4_{\kern 1.5pt\tilde{\kern -2pt i}=1}
\sum^4_{\tilde m=1}\sum^4_{\tilde n=1}\sum^4_{\tilde p=1}
\tilde\Omega[q]_{\kern 2pt\tilde{\kern -2pt i}\kern 0.2pt\tilde m\kern 0.1pt
\tilde n\tilde p}
\ \hat T^{\kern 2.5pt\tilde{\kern -2pt i}}_i
\ \hat T^{\kern 0.5pt\tilde m}_m\ \hat T^{\kern 0.5pt\tilde n}_n
\ \hat T^{\kern 0.5pt \tilde p}_p,\ q=2,\,\ldots,\,5.
\quad
\mytag{4.9}
\endgather
$$
This result is summarized in the following theorem.
\mytheorem{4.2} Under the right composition \mythetag{4.5} of a cubic 
transformation $\tilde f$ with a liner transformation $\varphi$ in 
\mythetag{1.2} its associated quartic forms $\omega[1]$, 
$\omega[2]$, $\omega[3]$, $\omega[4]$, $\omega[5]$, $\omega[6]$ are 
transformed according to the formulas \mythetag{4.7}, \mythetag{4.8}, 
and \mythetag{4.9}, where $\hat T$ is the block-diagonal matrix
similar to $\hat S$ in \mythetag{2.11} and built with the use of\/ $T$.
\endproclaim 
    The formulas \mythetag{4.7}, \mythetag{4.8}, \mythetag{4.9}
resemble the formula \mythetag{3.3} in Definition~\mythedefinition{3.1}.
For this reason the analogs of the forms $\omega[1]$, $\omega[2]$, 
$\omega[3]$, $\omega[4]$, $\omega[5]$, $\omega[6]$ in \mycite{14}
were called pseudotensors. However, this is not correct with respect 
to changes of coordinates of the form \mythetag{2.10} and 
\mythetag{2.11}.\par
    The origin and the prospective applications of the quartic forms
$\omega[1]$, $\omega[2]$, $\omega[3]$, $\omega[4]$, $\omega[5]$, $\omega[6]$ 
are associated with the following theorem.
\mytheorem{4.3} If a cubic transformation $f$ is produced as the
right composition \mythetag{4.5} of another cubic transformation 
$\tilde f$ with a liner transformation $\varphi$ of the form 
\mythetag{1.2}, then the determinants \mythetag{2.2} associated with
$f$ are expressed through the values of the quartic forms 
$\tilde\omega[1]$, $\tilde\omega[2]$, $\tilde\omega[3]$, $\tilde\omega[4]$, 
$\tilde\omega[5]$, $\tilde\omega[6]$ associated with the second cubic
transformation $\tilde f$ according to the formulas 
$$
\hskip -2em
\aligned
&G_{1111}(f)=\det T\cdot\tilde\omega[1](z^1,z^2),\\
&G_{1112}(f)=\det T\cdot\tilde\omega[2](z^1,z^2,z^3,z^4),\\
&G_{1122}(f)=\det T\cdot\tilde\omega[3](z^1,z^2,z^3,z^4),\\
&G_{1212}(f)=\det T\cdot\tilde\omega[4](z^1,z^2,z^3,z^4),\\
&G_{1222}(f)=\det T\cdot\tilde\omega[5](z^1,z^2,z^3,z^4),\\
&G_{2222}(f)=\det T\cdot\tilde\omega[6](z^3,z^4),
\endaligned
\mytag{4.10}
$$
where $z^1$, $z^2$, $z^3$, $z^4$ are given by the 
components of the matrix $T$ in \mythetag{1.2}:
$$
\xalignat 2
&\hskip -2em
z^1=T^1_1,
&&z^3=T^1_2,\\
\vspace{-1.5ex}
\mytag{4.11}\\
\vspace{-1.5ex}
&\hskip -2em
z^2=T^2_1,
&&z^6=T^2_2.\\
\endxalignat
$$
\endproclaim 
    Theorem~\mythetheorem{4.3} is proved by means of direct calculations 
which consist in deriving the formulas \mythetag{4.10} with the use of 
\mythetag{4.11}. 
\head
5. Conclusions.
\endhead
     Theorem~\mythetheorem{4.3} is the most important result of the 
present paper in view of its prospective application to classification 
of potentially invertible cubic transformations of the real plane 
$\Bbb R^2$. Such a classification of potentially invertible quadratic
transformations can be found in \mycite{14}.\par


\newpage
\Refs
\ref\myrefno{1}\by Vitushkin~A.~G.\paper On polynomial transformations
of\/ $\Bbb C^2$\inbook Proceedings of the International Conference on 
Manifolds and Related Topics in Topology, Tokyo, 1973\publ University
of Tokyo Press\yr 1975\pages 415--417
\endref
\ref\myrefno{2}\by Vitushkin~A.~G.\paper Some examples in connection 
with problems about polynomial transformations of\/ $\Bbb C^2$\jour
Izv. Akad. Nauk SSSR, Ser\. Mat.\vol 35\issue 2\yr 1971\pages 269--279
\endref
\ref\myrefno{3}\by Orevkov~S.~Yu.\paper On three-sheeted polynomial 
mappings of\/ $\Bbb C^2$\jour Izv\. Akad\. Nauk SSSR, Ser\. Mat.\vol 50
\issue 6\yr 1986\pages 1231--1240
\endref
\ref\myrefno{4}\by Domrina~A.~V., Orevkov~S.~Yu.\paper On four-sheeted 
polynomial mappings of\/ $\Bbb C^2$. I. The case of an irreducible ramification 
curve \jour Mat\. Zamretki\vol 64\issue 6\yr 1998\pages 847--862
\endref
\ref\myrefno{5}\by Guedj~V.\paper Dynamics of quadratic polynomial mappings 
of\/ $\Bbb C^2$\jour Michigan Math. Journal\vol 52\issue 3\yr 2004
\pages 627--648
\endref
\ref\myrefno{6}\paper
\myhref{http://en.wikipedia.org/wiki/Jacobian\podcherkivanie conjecture}
{Jacobian conjecture}\jour Wikipedia\publ Wikimedia Foundation Inc.\publaddr San 
Francisco, USA 
\endref
\ref\myrefno{7}\by Pinchuk~S\.~I\.\paper A counterexample to the strong real 
Jacobian conjecture\jour Mathematische Zeit\-schrift\vol 217\yr 1994\pages 1--4
\endref
\ref\myrefno{8}\by Jelonek~J\.~A\.\paper A geometry of polynomial transformations 
of the real plane\jour Bull\. Polish Academy of Sciences, Mathematics\vol 48
\issue 1\yr 2000\pages 57--62
\endref
\ref\myrefno{9}\by Campbell~L.~A.\paper On the rational real Jacobian 
conjecture \jour e-print \myhref{http://arxiv.org/abs/1210.0251}
{arXiv:1210.0251} in Electronic Archive \myEarXivlink
\endref
\ref\myrefno{10}\paper
\myhref{http://en.wikipedia.org/wiki/Euler\podcherkivanie 
brick}{Euler brick}\jour Wikipedia\publ 
Wikimedia Foundation Inc.\publaddr San Francisco, USA 
\endref
\ref\myrefno{11}\by Sharipov~R.~A.\paper Asymptotic estimates for roots 
of the cuboid characteristic equation in the nonlinear region\jour e-print
\myhref{http://arxiv.org/abs/1506.04705}{arXiv:1506.04705} in \myEarXivlink
\endref
\ref\myrefno{12}\by Sharipov~R.~A.\book Course of analytical
geometry\publ Bashkir State University\publaddr Ufa\yr 2011\moreref
see also \myhref{http://arxiv.org/abs/1111.6521}{arXiv:1111.6521}
in Electronic Archive \myEarXivlink
\endref
\ref\myrefno{13}\by Sharipov~R.~A.\book Quick introduction to tensor 
analysis\publ free on-line textbook 
\myhref{http://arxiv.org/abs/math/0403252}{arXiv:math/0403}
\myhref{http://arxiv.org/abs/math/0403252}{252}
in Electronic Archive \myEarXivlink
\endref
\ref\myrefno{14}\by Sharipov~R.~A.\paper A note on invertible quadratic 
transformations of the real plane\jour e-print 
\myhref{http://arxiv.org/abs/1507.01861}{arXiv:1507.01861} in Electronic 
Archive \myEarXivlink
\endref
\endRefs
\enddocument
\end